\documentclass[12pt]{article}
\usepackage{latexsym}
\usepackage{amssymb}
\usepackage{euscript}
\usepackage{amsmath}

\usepackage{graphicx}
\usepackage{color}

\def\mcc{M\raise.5ex\hbox{c}C}
\def\mccarthy{M\raise.5ex\hbox{c}Carthy}

\def\eg{{\it e.g. }}
\def\ie{{\it i.e. }}

\def\h{{\cal H}}

\def\M{{\cal M}}
\def\N{{\cal N}}

\def\MS{{\cal K}}

\def\ga{\gamma}

\def\a{\alpha}
\def\l{\lambda}
\def\z{\zeta}

\let\i=\infty

\def\la{\langle}
\def\ra{\rangle}
\def\={\ = \ }

\def\ba{{\cy7 B}}

\def\C{\mathbb C}
\def\R{\mathbb R}
\def\T{\mathbb T}
\def\bD{\mathbb D}
\def\B{\mathbb B}
\def\bH{\mathbb H}

\def\be{\setcounter{equation}{\value{theorem}} \begin{equation}}
\def\ee{\end{equation} \addtocounter{theorem}{1}}
\def\beq{\begin{eqnarray*}}
\def\eeq{\end{eqnarray*}}

\def\bp{{\sc Proof: }}
\def\ep{{}{\hfill $\Box$} \vskip 5pt \par}

\def\bl{\begin{lemma}}
\def\el{\end{lemma}}
\def\bt{\begin{theorem}}
\def\et{\end{theorem}}
\def\bprop{\begin{prop}}
\def\eprop{\end{prop}}
\def\bd{\begin{definition}}
\def\ed{\end{definition}}
\def\br{\begin{remark}}
\def\er{\end{remark}}
\def\bexer{\begin{exercise}}
\def\eexer{\end{exercise}}

\newtheorem{theorem}{Theorem}
\newtheorem{prop}[theorem]{Proposition}
\newtheorem{lemma}[theorem]{Lemma}
\newtheorem{definition}[theorem]{Definition}
\newtheorem{remark}[theorem]{Remark}
\newtheorem{cor}[theorem]{Corollary}

\newtheorem{question}[theorem]{Question}

\def\a{\alpha}

\def\N{\mathbb N}

\def\ba{\begin{eqnarray*}}
\def\ea{\end{eqnarray*}}

\def\mlt{\operatorname{Mult}}
\def\a{\mathfrak{a}}

\begin{document}
\title{Spaces of Dirichlet series with the complete Pick property}

\author{John E. M\raise.45ex\hbox{c}Carthy
\thanks{Partially supported by National Science Foundation Grant
DMS 1300280}
\and
Orr Moshe Shalit\footnote{Partially supported by ISF Grant 474/12 and
EU FP7/2007-2013 Grant 321749}}


\bibliographystyle{plain}

\maketitle

\begin{abstract}
We consider reproducing kernel Hilbert spaces of Dirichlet series with kernels of the form $k(s,u) = \sum a_n n^{-s-\bar u}$, and characterize when such a space is a complete Pick space.
We then discuss what it means for two reproducing kernel Hilbert spaces to be ``the same'', and introduce a notion of weak isomorphism.
Many of the spaces we consider turn out to be weakly isomorphic as reproducing kernel Hilbert spaces to the Drury-Arveson space $H^2_d$ in $d$ variables, where $d$ can be any number in $\{1,2,\ldots, \infty\}$, and in particular their multiplier algebras are unitarily equivalent to the multiplier algebra of $H^2_d$.
Thus, a family of multiplier algebras of Dirichlet series are exhibited with the property that every complete Pick algebra is a quotient of each member of this family.
Finally, we determine precisely when such a space of Dirichlet series is weakly isomorphic as a reproducing kernel Hilbert space to $H^2_d$ and when its multiplier algebra is isometrically isomorphic to $\mlt(H^2_d)$.
\end{abstract}

\section{Introduction}

Let $\h$ be a reproducing kernel Hilbert space on the set $X$, with kernel function $k(x,y)$ (below we shall also use the terminology {\em Hilbert function space on $X$}).
For any positive natural number $m$,
we say that $\h$ (or $k$) has the {\em $m$-Pick property} if, whenever $W_1, \dots, W_N$ is a finite set
of $m$-by-$m$ matrices, and $\l_1, \dots, \l_N$ are points in $X$, and the $mN$-by-$mN$ matrix given in block form by
\[
\Big[  k(\l_i, \l_j )\ [I_m - W_i W_j^* ]  \Big]
\]
is positive semi-definite, then there is a multiplier $\Phi$ in the closed unit ball of
${\rm Mult} (\h \otimes \C^m)$ that satisfies
\[
\Phi(\l_i) \= W_i, \qquad 1 \leq i \leq N .
\]
If $\h$ has the $m$-Pick property for all positive natural numbers, we say it has the {\em complete Pick property}.

The most well-known space with the complete Pick property is the Hardy space $H^2$, but there are others,
\eg \cite{ag90b, qui93, marsun, mcc04, hartz15}.
Spaces with the $m$-Pick property are described in \cite{agmc_loc}, but the description is cleaner for
spaces with the complete Pick property. These are totally described by the McCullough-Quiggin theorem
\cite{mccul92,qui93,agmc_cnp}. We say the kernel $k$ is {\em irreducible} if $X$ cannot be partitioned into two non-empty
sets $X_1, X_2$ such that $k(x,y) =0$ whenever $x \in X_1$ and $y \in X_2$.
We shall make a standing assumption throughout this note that all kernels are irreducible.

\bt
\label{thma1}
{\rm  (McCullough-Quiggin)} A necessary and sufficient condition for $k$ to have the complete Pick property
is that for any finite set $\{ \l_1, \dots, \l_N \}$ of distinct points in $X$, the matrix
\[
\left[
\frac{1}{k(\l_i, \l_j)}  \right]
\]
has exactly one positive eigenvalue.
\et
It was proved in \cite{agmc_cnp} that there is a universal
space with the complete Pick property, in the sense of Theorem~\ref{thma2} below.
For $d \in \{1,2,\ldots, \infty\}$, let $\B_d$ denote the open unit ball in a $d$-dimensional Hilbert space, and define a kernel $\a^d$ on $\B_d$
by
\[
\a^d(\zeta,\lambda) \= \frac{1}{1 - \langle \zeta, \lambda \rangle} .
\]
When $d = \infty$ we simplify notation to $\a^\infty = \a$.
Let $H^2_\infty$ denote the Hilbert function space on $\B_\infty$ with $\a$ as its reproducing kernel
(this is the infinite dimensional version of the Drury-Arveson space).
We denote by $\M_\infty$ the multiplier algebra $\mlt(H^2_\infty)$ of $H^2_\infty$.
The space $H^2_\infty$ not only has the complete Pick property (as can be easily seen from the McCullough-Quiggin theorem), but is universal.

We shall say that a kernel on $X$ is {\em normalized} at a point $\l_0 \in X$ if
$k(\zeta, \l_0) = 1$ for all $\zeta$. Any complete Pick kernel $k$  can be normalized by replacing it by
the equivalent kernel
\be
\label{eqa1}
\frac{ k(\zeta, \l)}{k(\z, \l_0) k(\l_0, \l)} ;
\ee
 the condition that $k$ is irreducible and has the $1$-Pick property means
that $k(\z,\l_0)$ is never $0$ \cite[Lemma 7.2]{ampi}, so \eqref{eqa1}
will remain holomorphic in $\z$ if $k$ is.

\bt
\label{thma2}
Suppose $k$ is a kernel normalized at $\l_0$. Then $k$ has the complete Pick property
if and only if there is a map $b: X \to \B_\infty$ that maps $\l_0$ to $0$ and satisfies
\[
k(\z,\l) \= \a (b(\z), b(\l)) .
\]
\et

It follows immediately from the theorem that every multiplier algebra of a complete Pick space is a quotient of $\M_\infty$.
The purpose of this note is to show that there is a space $\h$ of Dirichlet series that is also universal with
respect to having the complete Pick property, in the sense that every multiplier algebra of a complete Pick space is a quotient of $\mlt(\h)$. (This is Theorem~\ref{thm:universalPick} below).
On the other hand, the space $\h$ is not universal in the same sense as Theorem \ref{thma2} because its joint domain of definition is a half plane and this set turns out to be ``too small''.
Thus we will begin by exploring notions of isomorphism of reproducing kernel Hilbert spaces.

{\bf Acknowledgement:} We would like to thank Michael Hartz for pointing out a gap in an earlier version
of this paper.

\section{ When are two reproducing kernel spaces the same?}\label{sec:thesame}

A reproducing kernel Hilbert space $\h$ is equipped with a set $X$ on which the functions are defined.
But the functions are also defined on every subset of $X$, and perhaps there are sets containing $X$ to which the functions in $\h$ can be extended naturally.
In this section we discuss some issues arising from the fact that a reproducing kernel Hilbert space may be defined on a set which is not ``maximal''.
Our discussion is similar to the one carried out in \cite{cm95} (see their Definition 1.5, {\em algebraic consistency}) and \cite[Section 5]{hartz15}.

\begin{definition}
\label{defb1}
We say that the reproducing kernel Hilbert spaces $(\h_1, k_1)$ on $X_1$ and $(\h_2, k_2)$ on $X_2$
are {\em isomorphic as reproducing kernel Hilbert spaces} if there is a bijection $\phi : X_1 \to X_2$ such that
\[
k_2(\phi(x) , \phi(y)) = k_1(x,y) \quad \forall  x,y \in X_1 .
\]
\end{definition}
Equivalently, $(\h_1, k_1)$ on $X_1$ and $(\h_2, k_2)$ on $X_2$
are isomorphic as reproducing kernel Hilbert spaces if there is a bijection $\phi : X_1 \to X_2$ such that there is unitary isomorphism between the spaces that is induced by composition with the function $\phi$.
The above definition seems like the most natural condition for saying that two reproducing kernel Hilbert spaces are the same.
But the space $(\h_1, k_1)$ need not come presented with a maximal set $X_1$
on which it is a function space.
For example, suppose $X_1$ and $X_2$ are disjoint subsets
of the unit disk, both of which are sets of uniqueness, and let $\h_1$ be the Hardy space $H^2$
restricted to $X_1$, and $\h_2$ be $H^2$ restricted to $X_2$.
Both spaces are ``the same'', but the kernels seem to live on disjoint sets.

The points of $X_1$ are bounded point evaluations for $\h_1$. How does one find others?
If $x$ is a bounded point evaluation, then the kernel function $k_x$ is a joint eigenvector
for the adjoint of every multiplication operator; but this may not be the right generalization.
Consider for example the Fock space, all entire functions on $\C$ that are square integrable with respect
to the standard Gaussian measure. The only multipliers are the constant functions, so every vector
in the space is a joint eigenvector.

Instead, we shall use the following definition.
\begin{definition} A vector $v$ in the Hilbert function space $\h$
is a {\em generalized kernel function} if $v$ is non-zero and
\be
\label{eqa4}
\la fg, v \ra \= \la f, v \ra \, \la g, v \ra
\quad {\rm whenever\ } f, g \ {\rm and\ } fg \ {\rm are\ in\ }\h.
\ee
\end{definition}
Clearly, every kernel function is also a generalized kernel function.
Loosely, we think of the generalized kernel functions as being the evaluation functionals
for the largest set to which functions in $\h$ can be extended so that $\h$ continues to be a Hilbert function space.
The correct interpretation of the previous sentence requires caution: every Hilbert space can be considered as a Hilbert function space on itself -- a Hilbert space $H$ is just the set of all bounded linear functionals on $H$.
However, on this function space pointwise multiplication is not allowed (the product of two linear functionals is no longer linear), and one may say that it is not very interesting as a function space.

It should be kept in mind that when given a Hilbert function space $\h$ on a set $X$, the realization of elements of $\h$ as functions on $X$ determines a multiplication between elements of $\h$, and makes sense of the question whether $fg \in \h$ when $f,g \in \h$.
Thus the generalized kernel functions of $\h$ (given as a Hilbert function space on $X$) are the evaluation functionals for the largest set on which $\h$ is a Hilbert function space with the same algebraic structure determined by its realization as a function space on $X$.
Note also that if $\h$ is a Hilbert function space on $X$ and $f$ in $\h$ vanishes on $X$, then $f = 0 $ in $\h$.

In the most familiar examples of a Hilbert function space $\h$ on a set $X$, the only generalized kernel functions are the point evaluations (see Definition 1.5 and the surrounding discussion as well as Theorem 2.1.15 in \cite{cm95} and \cite[Section 5]{hartz15}).
We will encounter situations where $X$ is only a small part of the space of generalized kernel functions.

\bprop\label{prop:xhat}
Let $\h$ be a Hilbert function space on a set $X$ with kernel $k$.
Define a set $\hat{X} \subseteq \h$ by
\be\label{eq:xhat1}
\hat{X} = \{g : g \ \textrm{\em is  a generalized kernel function}\}.
\ee
Then there is a map $b: X \to \hat{X}$ (injective if $\h$ separates the points of $X$) and a Hilbert function space $\hat{\h}$ on $\hat{X}$, such that the map $\hat{f} \mapsto \hat{f} \circ b$ is an isometric isomorphism of $\hat{\h}$ onto $\h$.
Moreover, $\hat{\h}$ separates the points of $\hat{X}$, and the set of generalized kernel functions for $\hat{\h}$ is $\{ \hat{k} : \hat{k} \circ b \in \hat{X}\}$.
\eprop
\bp
Let $\hat{X}$ be as in \eqref{eq:xhat1}.
Define $b : X \to \hat{X}$ by
\[
b(x) = k_x.
\]
By definition $b$ is injective if $X$ separates points.
For every $f \in \h$, define $\hat{f} : \hat{X} \to \C$ by
\[
\hat{f}(g) = \langle f, g \rangle.
\]
We have that $\hat{f}(b(x)) = f(x)$.
Let $\hat{\h}$ be the space $\{\hat{f} : f \in \h\}$.
We see that $f \mapsto \hat{f}$ is a linear bijection that respects multiplication when defined, and its inverse is given by composition with $b$.
Defining the norm in $\hat{\h}$ by $\|\hat{f}\| := \|f\|$ makes these linear isomorphisms isometric.
Finally, it is clear that $\h$ and $\hat{\h}$ share the same set of generalized kernel functions, and those of $\h$ are identified with $\hat{X}$.
\ep
Note that $\h$ and $\hat{\h}$ might not be isomorphic as reproducing kernel Hilbert spaces (\ie in the sense of
Definition~\ref{defb1}), because $b$ may fail to be surjective.

On the other hand,  $b$ preserves some additional structure that $X$ has.
For example, if $X$ is a topological space and $k : X \times X \to \C$ is continuous then $b$ is continuous.
Indeed, this follows from considering
\[
\|k_x - k_y \|^2 = k(x,x) - k(x,y) - k(y,x) + k(y,y).
\]
Likewise, if $X \subseteq \C$ is a domain and $k : X \times X \to \C$ is continuous, holomorphic in the first variable and anti-holomorphic in the second variable, and if evaluation of the derivative at every point of $X$ is also bounded, then $b: X \to \h$ is anti-holomorphic.

In natural cases we can identify $\hat{X}$ with a concrete subset of $\C^d$.
\bprop
\label{prop9}
Let $\h$ be a Hilbert function space on a set $X$, and suppose that there are $d$ functions $\{\phi_1,
\phi_2,  \ldots \}$ (where $d \in \{1,2, \ldots, \infty\}$) such that the algebra generated by $\{\phi_1,
\phi_2,  \ldots \}$ is contained in $\h$ and dense in $\h$.
Then $\hat{X}$ can be identified with the set
\be\label{eq:xhat}
\hat{X} = \{(\langle \phi_1, g\rangle,  \langle \phi_2, g \rangle, \ldots) : g \ \textrm{\em is  a generalized kernel function}\}.
\ee
\eprop
\bp
This is similar to the previous proposition, with the change that we define $b : X \to \C^d$ by
\be\label{eq:b}
b(x) = (\phi_1(x),  \phi_2(x), \ldots) = (\langle \phi_1, k_x\rangle,  \langle \phi_2, k_x\rangle, \ldots) .
\ee
\ep

If one of the $\phi_i$'s is equal to the constant function $1$,
then it can be omitted in the above construction (provided the space is more than one-dimensional).
Thus if $X$ is a subset of $\C^d$ and $\h$ contains the algebra of polynomials as a dense subspace, then the set $\hat{X}$ of generalized kernel functions can also be identified with a subset of $\C^d$.

In many cases of interest, the function $b$ from \eqref{eq:b} can be chosen to play the role of the embedding $b \to \B_\infty$ from Theorem \ref{thma2} (see Section \ref{sec:uni_dir} below, or Section 7 in \cite{dhs15}).
Adapting the arguments of \cite[Section 5]{hartz15}
(which use slightly different definitions) we see that $\hat{X}$ can be identified with a subset of the smallest multiplier variety in the ball containing $b(X)$; we do not know whether in general $\hat{X}$ can be identified with the smallest multiplier variety containing $b(X)$.
See Remark \ref{rem:hat}.

\begin{definition}
Let $\h_1$  and $\h_2$ be  reproducing kernel Hilbert spaces.
A unitary $U : \h_1 \to \h_2$ is said to be a {\em multiplicative unitary} if
\be
\label{eqa3}
{\rm If\ } f,g\ \in \h_1, \ {\rm then\ } fg \ \in \h_1 \ {\rm if \ and\ only\ if \ }
U(f) U(g) \ \in \h_2 .
\ee
\be
\label{eqa2}
U(fg) \= U(f) U(g) \quad {\rm whenever\ } f, g \ {\rm and\ } fg \ {\rm are\ in\ }\h_1.
\ee
\end{definition}
Note that $U$ is a multiplicative unitary if and only if $U^*$ is.
\begin{definition}
Let $\h_1$  and $\h_2$ be  reproducing kernel Hilbert spaces.
We shall say that $\h_1$ and $\h_2$ are {\em weakly isomorphic as reproducing kernel Hilbert spaces}
if there is a multiplicative unitary $U $ from $\h_1 $ onto $\h_2$.
\end{definition}
\bprop\label{prop:isotohat}
For every reproducing kernel Hilbert space $\h$ on a set $X$, $\h$ is weakly isomorphic as a reproducing kernel Hilbert space to $\hat{\h}$ on $\hat{X}$.
\eprop
This follows from Proposition~\ref{prop:xhat}.

We shall prove (Theorem~\ref{thm:universalPick}) that there are Hilbert spaces of Dirichlet series that are weakly isomorphic
as reproducing kernel Hilbert spaces to $H^2_\infty$.

\bprop
Let $U : \h_1 \to \h_2$ be a unitary.
Then $U$ is multiplicative if and only if it maps the set of generalized kernel functions in $\h_1$ onto
the set of generalized kernel functions in $\h_2$.
When these conditions hold, then there is a multiplicative unitary $\hat{U} : \hat{\h}_1 \to \hat{\h}_2$ and a bijection $\phi : \hat{X}_1 \to \hat{X}_2$ such that $\hat{U}(f) = f \circ \phi^{-1}$.
\eprop
\bp
Suppose that $U$ is a multiplicative unitary.
Let $v$ be a generalized kernel function in $\h_1$, and let $V = U v$.
Let $F, G $ and $FG$ be in $\h_2$. Let $f = U^* F$ and $g = U^* G$.
Then $fg \in \h_1$, by \eqref{eqa3}.
We have
\beq
\la FG , V \ra_{\h_2} &\= &
\la U^*(FG) , U^* V \ra_{\h_1} \\
&=& \la fg, v \ra \\
&=& \la f, v \ra \ \la g, v \ra \\
&=&\la F, V \ra_{\h_2} \ \la G, V \ra_{\h_2}.
\eeq
Thus $U$ maps $\hat{X}_1$ into $\hat{X}_2$; applying this to $U^*$ we conclude ``onto''.

Conversely, assume that $U$ maps $\hat{X}_1$ onto $\hat{X}_2$.
Applying Proposition \ref{prop:isotohat} we pass to a unitary $\hat{U}: \hat{\h_1}
\to \hat{\h}_2$ which maps $\hat{X}_1$ onto $\hat{X}_2$.
There exists therefore a bijective function $\phi : \hat{X}_1 \to \hat{X}_2$ such that
$\hat{U} k^1_{x} = k^2_{\phi(x)}$ for every $x \in \hat{X}_1$.
It follows that $\hat{U}^*$ is implemented by composition with $\phi^{-1}$, therefore $\hat{U}$ is implemented by composition with $\phi$.
In particular, $\hat{U}$ is multiplicative, so $U$ is too.
\ep

\begin{cor}
Let $\h_1$ and $\h_2$ be reproducing kernel Hilbert spaces.
Then $\h_1$ and $\h_2$ are weakly isomorphic as reproducing kernel Hilbert spaces if and only if $\widehat{\h}_1$ and $\widehat{\h}_2$ are isomorphic as reproducing kernel Hilbert spaces.
\end{cor}

\section{The complete Pick property for spaces with nice bases}

\bd Let $Y$ be a set. A sequence $\{\phi_n\}_{n=1}^\infty$ of functions defined on $Y$  is said to be {\em strongly linearly independent on $Y$} if for all $k \in \N$, there is no series of the form $\sum_{n\neq k} c_n \phi_n$ which  converges pointwise to $\phi_k$ on $Y$.
\ed
Examples of strongly independent sequences are given by $\phi_n(z) = z^n$ or $\phi_n(s) = n^{-s}$
on non-empty open sets. In fact, any space of functions where there is a series expansion with a uniqueness theorem would be an example.

\bl
Let $\{\phi_n\}_{n=1}^\infty$ be a strongly linearly independent sequence of functions on $Y$, and let $\{a_n\}_{n=1}^\infty$ be a sequence of real numbers. Consider the kernel
\[
K(x,y) = \sum_{n=1}^\infty a_n \phi_n(x) \overline{\phi_n(y)}.
\]
Then $K$ is positive semi-definite if and only if $a_n \geq 0$ for all $n$.
\el
\bp
We prove the nontrivial (but probably well known) direction: suppose that $K$ is positive semi-definite.
Put $I = \{n : a_n <0\}$ and $J = \{n : a_n \geq 0\}$. Then
\[
K_I(x,y) := \sum_{n\in I}-a_n  \phi_n(x) \overline{\phi_n(y)}
\]
and
\[
K_J(x,y) := \sum_{n \in J} a_n \phi_n(x) \overline{\phi_n(y)}
\]
are both positive semi-definite kernels, and we have the kernel inequality $K_I \leq K_J$.
In particular, for every $m \in I$
\[
\phi_m(x) \overline{\phi_m(y)} \leq cK_J(x,y)
\]
for some positive constant $c$. We deduce (by \cite[Theorem 4.15]{Pau09notes}) that $\phi_m$ is in the reproducing kernel Hilbert space $H(K_J)$ determined by $K_J$. By \cite[Theorem 3.12]{Pau09notes} it follows that $\{\phi_n\}_{n \in J}$ is a Parseval frame for $H(K_J)$. From \cite[Proposition 3.10]{Pau09notes} it then follows that
\[
\phi_m = \sum_{n \in J} \langle \phi_m, \phi_n \rangle \phi_n
\]
in norm (hence pointwise), contradicting the assumption that the sequence $\{\phi_n\}_{n=1}^\infty$ is strongly linearly independent. If follows that $I$ must be empty.
\ep

Let $\h$ be a reproducing kernel Hilbert space on a set $X$ that contains the constants. Let $\{\phi_n\}_{n=1}^\infty$ be an orthogonal basis for $\h$.
Then the kernel of $\h$ is given by
\be\label{eq:kernelphi}
k(x,y) = \sum_{n=1}^\infty a_n \phi_n(x) \overline{\phi_n(y)},
\ee
where $a_n = \|\phi_n\|^{-2}$.

\bprop\label{prop_pick} Suppose that $k$ is a kernel on $X$ that is never zero, and is normalized
at $x_0$.
Suppose that  $\{\phi_n\}_{n=1}^\infty$ is an orthogonal basis for $\h$,
and that the sequence $\{\phi_n\}_{n=1}^\infty$ is such that one can write $1 - k(x,y)^{-1}$ as
\be\label{eq:kernelinv}
1 - k(x,y)^{-1} = \sum_{n=1}^\infty \alpha_n \phi_n(x) \overline{\phi_n(y)}.
\ee
Then $k$ is a complete Pick kernel if $\alpha_n \geq 0$ for all $n \geq 1$. Conversely, if $\{\phi_n\}_{n=1}^\infty$ is a strongly linearly independent sequence on $X$, then the condition $\alpha_n \geq 0$ is also necessary for that $\h$ be a complete Pick space.
\eprop
\bp
Suppose that $\alpha_n \geq 0$ for all $n \geq 1$. Put $b_n = \sqrt{\alpha_n}$. Define a function
$f : X \rightarrow \B_\infty$ by
\[
f(x) = (b_1 \phi_1(x), b_2 \phi_2(x), \ldots).
\]
By (\ref{eq:kernelinv}) and positivity of $k$, we have $\|f(x)\|^2 = 1 - k(x,x)^{-1} <1$ for all $x$, so $f$ indeed maps
into $\B_\infty$.
Rearranging (\ref{eq:kernelinv}), we find
\[
k(x,y) = \frac{1}{1- \langle f(x), f(y) \rangle}.
\]
It follows that $\h$ is isomorphic to the space $\overline{\textrm{span}}\{\a_{f(x)} : x \in X\}$, thus $\h$ is a complete Pick space, since $H^2_\infty$ is.

Conversely, assume that $\h$ is a complete Pick space and that $\{\phi_n\}_{n=1}^\infty$ is strongly linearly independent.
We re-organize (\ref{eq:kernelinv}) as
\be\label{eq:koneminussum}
k(x,y)^{-1} = 1 - \sum_{n=1}^\infty \alpha_n \phi_n(x) \overline{\phi_n(y)}.
\ee
By the McCullough-Quiggin Theorem~\ref{thma1}, the complete Pick property implies that for every choice of points $x_1, \ldots, x_N \in X$, the matrix $\left[k^{-1}(x_i,x_j)\right]_{i,j=0}^N$ has exactly one positive eigenvalue.
By taking the Schur complement with respect to the $(0,0)$ entry, we
conclude that the $N$-by-$N$ matrix with entries
\[
\left[ \frac{1}{k(x_i, x_j)} \ - \ \frac{k(x_0, x_0)}{k(x_i, x_0) k(x_0, x_j)} \right]_{i,j=1}^N
= \left[ \frac{1}{k(x_i, x_j)} \ - \ 1 \right]_{i,j=1}^N
\]
is negative semi-definite. Using \eqref{eq:koneminussum}, we get
\[
\left[\sum_{n=1}^\infty \alpha_n \phi_n(x_i) \overline{\phi_n(x_j)}\right]_{i,j=1}^N
\ = \
\left[ 1 \ - \ \frac{1}{k(x_i, x_j)} \right]_{i,j=1}^N \ \geq \ 0.
\]
Thus the kernel
\[
\tilde{k}(x,y) := \sum_{n=1}^\infty \alpha_n \phi_n(x) \overline{\phi_n(y)}
\]
is positive semi-definite on $X$. By the lemma, $\alpha_n \geq 0$ for all $n$.
\ep

\section{Spaces of Dirichlet series}\label{subsec:spaces_Dirichlet}

We can apply  Proposition \ref{prop_pick} to spaces of Dirichlet series.
Let us provide some details.

Let $\h$ be a Hilbert function space of Dirichlet series,
with the kernel given by $k(s,u) = \sum a_n n^{-s-\bar u}$, and suppose that this kernel converges for all $s,u$ in some half space $\bH_\delta := \{ s : \Re(s) > \delta\}$.
For simplicity we assume that $a_1 = 1$.
Since the coefficients of this series are positive, the abscissae of absolute and uniform convergence are
the same as the abscissa of convergence of the series.

One sees that $\h$ is the space of all functions $h$ with Dirichlet series $h(s) = \sum \ga_n n^{-s}$ satisfying $\sum a_n^{-1} |\ga_n|^2 < \infty$ (the formal Dirichlet series $k_u(s) = \sum a_n n^{-\bar u} n^{-s}$ is readily seen to satisfy this bound provided $\textup{Re}\, u > \delta$).
Consequently, the Dirichlet series of every $h \in \h$ converges also in $\bH_\delta$.
It follows that for every Hilbert space of Dirichlet series $\h$ as above there exists a $\delta_0 \in [-\infty, \infty)$ such that the Dirichlet series for every $h \in \h$ converges in the half plane $\bH_{\delta_0}$, and that there are functions in $\h$ that do not converge on any strictly bigger half plane.

Let $c_n$ denote the coefficients of $k^{-1}$ as follows:
\[
\frac{1}{\sum_{n=1}^\i a_n n^{-s}} \= \sum_{n=1}^\i c_n n^{-s} .
\]
One finds that $c_1 = 1$ and that $c_n$ for $n>1$ are given by the recursive formula:
\be\label{eq:cn}
c_n = - \sum_{d<n, d|n} a_{n/d} c_d .
\ee
Now we will see how to obtain that $\h$ is a complete Pick space if and only if $c_n \leq 0$ for all $n > 1$.

Suppose that $c_n \leq 0$ for all $n>1$. Then the sum $\sum_{d<n, d|n} a_{n/d} c_d$ is non-negative for all $n>1$. Since the first term is equal to $a_n$ and $a_n \geq 0$, and as all other terms are negative, it follows that $|\sum_{d<n, d|n} a_{n/d} c_d| \leq  a_n$ for all $n>1$. Hence $|c_n| \leq a_n$ for all $n$.

We conclude the following: if $c_n \leq 0$ for $n>1$ then $|c_n| \leq a_n$ for all $n$; in particular, the Dirichlet series for $k^{-1}$ converges  on $\bH_{\delta}$ too, and Proposition \ref{prop_pick} applies to show that $\h$ is a complete Pick space.

Conversely, if $\h$ is a complete Pick space with kernel $k(s,u) = f(s+\bar u)$ where $f(s) = 1 + \sum_{n \geq 2} a_n n^{-s}$ converges on some half plane $\bH_{\sigma_1}$, then $k^{-1}$ is also given by a Dirichlet series that converges uniformly on some half plane, say $\bH_{\sigma_2}$.
Now put $\sigma = \max\{\sigma_1, \sigma_2\}$. It is well known that the Dirichlet series of a function is unique, hence the condition of strong linear independence is satisfied.
Now Proposition \ref{prop_pick} applies to $\h\big|_{\bH_\sigma}$ (which is still a complete Pick space), and we conclude that $c_n \leq 0$ for $n>1$ (and the remarks above now show that $\sigma = \sigma_1$).
Thus, we have proved Theorem \ref{thm1}.

\bt
\label{thm1}
Suppose $\h$ is a holomorphic Hilbert space  with kernel
function
$$
k(s,u) \= \sum_{n=1}^\i a_n n^{-(s+\bar u)}
$$
and assume $a_1 \neq 0$.
Let the Dirichlet coefficients
of $\frac{1}{k}$ at infinity be given by
\be
\label{eqgb3}
\frac{1}{\sum_{n=1}^\i a_n n^{-s}} \= \sum_{n=1}^\i c_n n^{-s} .
\ee
Then $\h$ has the complete Pick property if and only if
$$
c_n \leq 0 \qquad \forall \, n \geq 2 .
$$
\et

Examples of kernels with the complete Pick property are easy to come by using Theorem \ref{thm1} and known formulas for
the zeta function \cite{tit86}. For example, let $k(s,u) = \phi( s + \bar u)$. Then this will give a complete Pick kernel if
\beq
\phi(s) &=& \frac{1}{2 - \zeta (s)} = \sum_{n=1}^\i \frac{f(n)}{n^s}, \\
\phi(s) &=& \frac{\zeta(s)}{\zeta(s) + \zeta'(s)}, \\
\phi(s) &=& \frac{\zeta(2s)}{2 \zeta(2s) - \zeta(s) }.
\eeq
In the first formula, $f(n)$ is the number of distinct ways in which $n$ can be factored  (where the order matters).

We finish this section by showing that Hilbert spaces of Dirichlet series with the complete Pick property cannot be supported on the entire plane.
\bt
Let $\h$ be a Hilbert space of Dirichlet series with kernel $k(s,u) = 1+\sum_{n\geq 2} a_n n^{-s-\bar u}$, and suppose that $\h$ has the complete Pick property, and that $\dim \h > 1$.
If $\bH_{\delta}$ is the largest half plane of convergence for $\h$, then $\delta > -\infty$.
\et
\bp
If $\delta = -\infty$, then for every $t >0$ the series $\sum a_n n^{t}$ converges absolutely.
Thus for all $t > 0$ there exists a constant $M_t$ such that $a_n \leq M_t n^{-t}$ for all $n\geq 1$.
On the other hand, inverting \eqref{eq:cn} and using Theorem \ref{thm1}, we find that
\[
a_n = \sum_{d<n, d|n}a_d |c_{n/d}|,
\]
and therefore
\[
a_{n^k} \geq |c_n|^k.
\]
Finding some $n>1$ so that $c_n \neq 0$, we find
\[
|c_n|^k \leq a_{n^k} \leq M_t n^{-kt}
\]
for all $k$ and all $t$.
Fixing $n$ and taking $t$ sufficiently large we obtain a contradiction.
\ep

It follows from the above theorem and a change of variables that all reproducing kernel Hilbert spaces of Dirichlet series with the complete Pick property are isomorphic as reproducing kernel Hilbert spaces to a space with joint domain of convergence equal to the right half plane $\bH_0$.

\section{The universal representation}\label{sec:uni_dir}

Let $\{b_k\}_{k=1}^\i$ be a sequence of positive numbers such that $\sum_{k=1}^\i b^2_k = 1$.
Motivated by the first part of the proof of Proposition \ref{prop_pick}, we consider the map
$f : \bH_0 \rightarrow \B_\i$ given by
\be\label{eq:b_embed}
f(s) = (b_1 p_1^{-s}, b_2 p_2^{-s}, b_3 p_3^{-s}, \ldots),
\ee
where $p_k$ denotes the $k$th prime.
We define a kernel in $\bH_0$ by
\[
k(s,u) = \a(f(s),f(u)) = \sum_n a_n n^{-s-\bar u} ,
\]
and we denote by $\h$ the Hilbert function space determined by $k$.
We have that $\h$ is a complete Pick space on $\bH_0$.

Let us recall some familiar facts about $\h$ (see \cite{salsha15}).
The space $\h$ is isometric to the restriction of $H^2_\infty$ to the smallest multiplier-variety $V$ in $\B_\infty$ that contains $f(\bH_0)$, which means
\be\label{eq:hull}
V = \{z \in \B_\infty : g(z) = 0 \textrm{ for all } g \in \M_\infty \textrm{ such that } g\big|_{f(\bH_0)}\equiv 0\}.
\ee
The mapping $U: k_\lambda \mapsto \a_{f(\lambda)}$ extends to a unitary map from $\h$ onto $\MS_{f(\bH_0)}$, where
\[
\MS_{f(\bH_0)} := \vee \{\a_{f(s)} : s \in \bH_0\} \subseteq H^2_\infty.
\]

Denoting 
\[
\MS_V := \vee \{\a_v : v \in V\} \subseteq H^2_\infty, 
\]
we have that, as subspaces of $H^2_\infty$, $\MS_V = \MS_{f(\bH_0)}$ (see Proposition 2.2 in \cite{drs15}), 
but it is important to remember that we consider $\MS_{f(\bH_0)}$ as Hilbert function space on $f(\bH_0)$ and $\MS_V$ as a Hilbert function space on $V$.

The adjoint of $U$ is given by $U^*h = h \circ f$.
It is clear that $U$ is a multiplicative unitary from $\h$ onto $\MS_{f(\bH_0)}$;
we will show below (Theorem~\ref{thm:universalPick})  that $U$ is a multiplicative unitary from $\h$ onto $\MS_V$.
Note the difference: we must show that if $g_1,g_2$ and $h$ are in $H^2_\infty$, and $g_1 g_2 = h$ on $f(\bH_0)$,
then $g_1 g_2 = h$ on $V$.

%

Let $\mlt(\h)$ be the multiplier algebra of $\h$.
Then $\mlt(\h)$ is isomorphic to $\M_V := \mlt (\MS_V)$; the  isomorphism $\Phi : \M_V \rightarrow \mlt(\h)$ is given by
\[
\Phi (M_g) = U^* M_g U, \quad  \quad g \in \M_V.
\]
For $g \in \M_V$ and $h \in \MS_V$, we compute
\[
U^* (g h) = (g \circ f) \cdot (h \circ f) = (g \circ f) U^* h.
\]
Thus $\Phi(M_g) = U^* M_g U = M_{g \circ f}$, or simply $\Phi(g) = g \circ f$.

It is interesting to see where $\Phi^{-1}$ sends the functions $n^{-s}$.
To this end, we compute
\[
\Phi (z_k)(s) = z_k \circ f(s) = b_k p_k^{-s}.
\]
Thus $\Phi^{-1}(p_k^{-s}) = b_k^{-1} z_k$, and $\Phi^{-1}(n^{-s})$ is given by the appropriate product, determined by the prime factoring of $n$.
To set notation we spell this out: if $n = p_1^{\mu_1} \cdots p_k^{\mu_k}$, we write $\mu(n) = \mu = (\mu_1, \ldots, \mu_k,0,0,\ldots)$, we write $ n(\mu) = n$, and we have
\[
\Phi^{-1}(n^{-s}) = (b^{\mu(n)})^{-1}  z^{\mu(n)}.
\]

\bt
\label{thm:universalPick}
Let the notation be as above.
Then $V = \B_\infty$, and the map $U: k_\lambda \mapsto \a_{f(\lambda)}$ extends to a multiplicative unitary from $\h$ onto $\MS_V = H^2_\infty$.
Thus, $\h$ is weakly isomorphic as a reproducing kernel Hilbert space to $H^2_\infty$, and in particular $\mlt(\h)$ is unitarily equivalent to $\M_\infty :=\mlt(H^2_\infty)$.
\et
We will need two lemmata.
For $r \in (0,1)$, and $g : \B_\infty \to \C$, we denote by $g_r$ the function $g_r(z) = g(rz)$ for all $z \in \B_\infty$.

\bl\label{lem:ineqball}
Let $g : \B_\infty \to \C$ be a function such that $g_r \in \M_\infty$ for all $r \in (0,1)$.
For every $\rho \in (0,1)$, there exists a constant $C_\rho$ such that for all $z,w \in \B_\infty$,
\be\label{eq:ineqball}
\|z\|,\|w\| < \rho \Rightarrow  \|g(w) - g(z)\|\leq C_{\rho} \|z - w\|.
\ee
\el

\bp
We begin by proving the result for $g \in \M_\infty$.
Without loss of generality, we may assume that $\| g \|_{\mlt(H^2_d)} \leq 1$.
By the positivity of the 2 point Pick matrix of $g$, we have
\begin{eqnarray*}
\frac{(1 - |g(z)|^2)(1- |g(w)|^2)}{(1 - \| z \|^2)(1- \| w \|^2)}
&\geq &
\frac{| 1 - g(z) \overline{g(w)} |^2}{| 1 - \langle z, w \rangle |^2} \\
\therefore
\left| \frac{g(z) -g(w)}{1 - g(z) \overline{g(w)}} \right|^2
& \leq &
1 - \frac{(1 - \| z \|^2)(1- \| w \|^2)}{| 1 - \langle z, w \rangle |^2}
\\
&\leq&
\left\|
\frac{ z  - w }{1 - \langle z, w \rangle} \right\|^2 .
\end{eqnarray*}
Since $\|1 - \langle z, w \rangle\|  \geq 1 - \rho^2$ and $|1 - g(z) \overline{g(w)}| \leq 2$, we obtain the result for multipliers.

Now let $g$ be as in the statement of the lemma.
Fix $r \in (\rho , 1)$, and consider $g_r$.
Then $g_r$ is a multiplier, and by the previous paragraph there is a constant $C'$ such that for al $z,w \in \B_\infty$,
\[
\|z\|,\|w\| < \rho/r \Rightarrow  \|g(rw) - g(rz)\|\leq C' \|z - w\|.
\]
Setting $C_\rho = C'/r$, we obtain \eqref{eq:ineqball}.
\ep

Note that the hypotheses of Lemma~\ref{lem:ineqball}  hold for all $g \in H^2_\infty$,
and therefore for all $g \in \M_\infty$,
because $g_r$ then has a Taylor series that converges absolutely in a neighborhood of the ball for all $r$.

\bl\label{lem:vanish}
Let $g : \B_\infty \to \C$ be a function such that $g_r \in \M_\infty$ for all $r \in (0,1)$.
If $g$ vanishes on $f(\bH_0)$, then $g = 0$ on all $\B_\infty$.
\el
\bp
Fix $\epsilon > 0$ and define $L_\epsilon$ to be the line
$L_\epsilon = \{\epsilon + it : t \in \R\}$.
For every $N$ we let $P_N$ denote the orthogonal projection onto $\{z \in \ell^2 : z_{N+1} = z_{N+2} = \ldots = 0\}$.
We define $g_N$ to be the restriction of $g$ to the subspace $P_N \B_\infty = \{z \in \B_\infty : z_{N+1} = z_{N+2} = \ldots = 0\}$.  It will be convenient to let $W $ denote $ f(L_\epsilon)$, and let $W_N = P_N f(L_\epsilon)$.
Thus
\[
W_N = \{(b_1 p_1^{-\epsilon} e^{-i \log p_1 t} , \ldots, b_N p_N^{-\epsilon} e^{-i \log p_N t}) : t \in \R\} .
\]
Suppose  we  can show that $g_N$ vanishes on $W_N$ for all $N$.
By Kronecker's theorem, $W_N$
is dense in the polytorus $b_1 p_1^{-\epsilon} \T \times \cdots \times b_N p_N^{-\epsilon} \T \subset \B_N$.
So if $g_N$ vanished on $W_N$, by the maximum principle it would vanish on a polydisk,
and hence it would be zero on the whole ball $\B_N$.
 But $g$ has a power series of the form
\[
g(z) = \sum c_\mu z^\mu,
\]
where $\mu$ runs over all finite multi-indices.
From $g_N \equiv 0$ for all $N$, it would follow that $c_\mu = 0$ for all $\mu$, hence $g \equiv 0$.

Thus we must show that $g_N$ vanishes on $W_N$ for all $N$.
Fix $N_0$, and let $N > N_0$.
Now apply Lemma \ref{lem:ineqball} to $g$ with $w \in W$ and $z = P_N w$, noting that
for any $ 2^{-\epsilon} < \rho < 1$, we have
 $W \subset \rho \B_\infty$.
Since $g(w) = 0$ and $g(z) = g_N(z)$, we get
\[
| g_N(z) |  \leq  C_\rho \|w - z\| = C_\rho \sqrt{\sum_{k=N+1}^\i |b_k|^2 p_k^{-2\epsilon}}
.
\]
This gives that $|g_N|$ restricted to $W_N$ has values less than $C_\rho r_N$, where $$
r_N \ := \ \sqrt{\sum_{k>N}b_k^2 p_k^{-2\epsilon}} .$$
Since $W_N$ is dense in the polytorus $$b_1 p_1^{-\epsilon} \T \times \cdots \times b_N p_N^{-\epsilon} \T, $$
 the maximum principle gives that $\|g_N\|_\infty \leq C_\rho r_N$ on the polydisk
 $$ b_1 p_1^{-\epsilon}\overline{ \bD} \times \cdots \times  b_N p_N^{-\epsilon} \overline{\bD} .
 $$
Since $N_0 < N$, it follows that
\be
\label{eqd1}
|g_{N_0}| \ \leq\  C_\rho r_N
\ee
 on $W_{N_0}$. Letting $N \to \i$ in \eqref{eqd1}, we get $g_{N_0} \equiv 0$ on $W_{N_0}$, as required.
\ep

\noindent{\sc Proof of Theorem \ref{thm:universalPick}: }
To show that  the multiplier closure $V$ of $f(\bH_0)$ (given by \eqref{eq:hull} ) is equal to $\B_\infty$,
we must 
show that the only multiplier $g$ that vanishes on $f(\bH_0)$ is $g = 0$. This follows from Lemma~\ref{lem:vanish}.

It remains to show that $U$ is a multiplicative unitary.
Clearly, $U^*$ is a multiplicative unitary from $\MS_{f(\bH_0)}$ to $\h$, because it is implemented by composition with a bijective function $f : \bH_0 \to f(\bH_0)$.
Now, every $h \in \MS_{f(\bH_0)}$ (considered as reproducing kernel Hilbert space on $f(\bH_0)$) extends uniquely to a function in $\MS_V = H^2_\infty$ (a reproducing kernel Hilbert space on $\B_\infty$).
We need to show that this extension operator is a multiplicative unitary.

When viewing $\MS_{f(\bH_0)}$ as the  subspace of $H^2_\infty$ spanned by $\a_{f(s)}$ ($s \in \bH_0$), this extension operator becomes the identity map, so it may seem like there is nothing to prove.
But there is something to prove:
we have to show that if $g_1, g_2,h \in H^2_\infty$, and that $g_1 g_2 \big|_{f(\bH_0)} = h\big|_{f(\bH_0)}$, then $g_1 g_2 = h$ on the whole ball.
This would show that $g_1 g_2 \in H^2_\infty$ if and only if $g_1 \big|_{f(\bH_0)} g_2 \big|_{f(\bH_0)} \in \MS_{f(\bH_0)}$.

The function $F = g_1 g_2 - h$ satisfies the assumptions of Lemma \ref{lem:vanish} (since $h_r$, $(g_1)_r$ and $(g_2)_r$  are multipliers, $F_r = (g_1)_r (g_2)_r - (h)_r$ is a multiplier too), so invoking this lemma
completes the proof.
\ep
\begin{remark}\label{rem:hat}
\emph{It follows from the proof of \cite[Lemma 5.2(a)]{hartz15} (where the setting is somewhat different),
or from Proposition~\ref{prop9},
 that the generalized kernel functions of $H^2_\infty$ are precisely the point evaluations in $\B_\infty$, thus in the setting of Theorem \ref{thm:universalPick} we have the identification $\widehat{f(\bH_0)} = V = \B_\infty$.}
\end{remark}

\begin{cor}
\label{cord1}
The norm of an element $h(s) = \sum \gamma_n n^{-s} \in \h$ is given by
\be
\label{eqd2}
\|h\|^2 = \sum_n \frac{|\gamma_n|^2}{(b^{\mu(n)})^{2}} \frac{\mu(n)!}{|\mu(n)|!}.
\ee
\end{cor}
\bp
We have
\beq U^* (z^\mu) (s) & \ = \ & z^\mu \circ f(s) \\
& =& b^\mu n(\mu)^{-s},
\eeq
so $$
U(n^{-s}) \ = \ \frac{1}{b^{\mu(n)}} z^{\mu(n)} .
$$
Using $\|z^{\mu}\|^2 = \frac{\mu(n)!}{|\mu(n)|!}$ we get \eqref{eqd2}.
\ep
Note that comparison of \eqref{eqd2} with $\|h\|^2 = \sum |\ga_n|^2 a_n^{-1}$ yields the formula
\[
a_n = b^{2\mu(n)} \frac{|\mu(n)|!}{\mu(n)!},
\]
in agreement with the inversion formula for Dirichlet series.

Using the universal property of the shift on Drury-Arveson space \cite{dru78}, we also obtain the following von Neumann type inequality.
\begin{cor}
Let $T = (T_1, \ldots, T_d)$ be a commuting row contraction.
Then for every polynomial $Q$ in $d$ variables, one has
\begin{equation}\label{eq:vN_inequality}
\|Q(T_1, \ldots, T_d)\| \leq \|Q(b_1 p_1^{-s}, \ldots, b_d p_d^{-s})\|_{\mlt(\h)} .
\end{equation}
\end{cor}

We get a universal kernel from any choice of positive sequence $b_k$ that satisfies
$\sum b_k^2 = 1$.
Here is a particular choice that gives a nice form for the kernel function.
Let
\[
P(s) \= \sum_{p \ \rm prime} p^{-s}
\]
be the {\em prime zeta function}.
Let
\[
b_k \= \frac{1}{ \sqrt{P(2)} p_k} .
\]
Then by Corollary \ref{cord1} we get that
\beq
k(s,u) &\= &
\frac{P(2)}{P(2) - P(2 + s + \bar u)} \\
& =  & \a(f(s),f(u)),
\eeq
where $f(s) = (b_1 2^{-s}, b_2 3^{-s}, b_3 5^{-s}, \ldots)$,
is a universal complete Pick kernel, in the sense that every complete Pick space is a quotient of $\h(k)$ and every complete Pick algebra is the quotient of $\mlt(\h(k))$.
However, it is not universal in the sense of Theorem \ref{thma2}, and in particular $\h$ is not isomorphic as a reproducing kernel Hilbert space to $H^2_\infty$.

\section{Which complete Pick spaces of Dirichlet series are universal?}\label{sec:which}

In Section \ref{sec:uni_dir} we saw that some particular spaces of Dirichlet series are weakly isomorphic as reproducing kernel Hilbert spaces to $H^2_\infty$, and that their multiplier algebras are unitarily equivalent to $\M_\infty$.
In the present section we ask which of the complete Pick spaces of Dirichlet series have this property.

Fix $d \in \{1,2,\ldots, \infty\}$, and let $\{b_k\}_{k=1}^d$ be a sequence of positive numbers such that $\sum_{k=1}^d b^2_k = 1$.
Let $n_1, n_2, \ldots, $ be an increasing sequence of positive integers and define the map
$f : \bH_0 \rightarrow \B_d$ by
\be\label{eq:b_embedn}
f(s) = (b_1 n_1^{-s}, b_2 n_2^{-s}, b_3 n_3^{-s}, \ldots) .
\ee
Letting $\a^d$ denote the kernel of the space $H^2_d$, we define a kernel in $\bH_0$ by
\[
k(s,u) = \a^d(f(s),f(u)) = \frac{1}{1 - \la f(s), f(u) \ra} = \sum_n a_n n^{-s-\bar u} ,
\]
and we denote by $\h$ the Hilbert function space determined by $k$.
By Proposition \ref{prop_pick}, we have that $\h$ is a complete Pick space on $\bH_0$.

\bt
\label{thmm3}
The multiplier closure of
$f( \bH_0)$ is equal to $\B_d$ if and only if the sequence $\log n_1, \log n_2, \ldots$ is linearly independent over $\mathbb Q$.
In this case $\h$ is weakly isomorphic as a reproducing kernel Hilbert space to $H^2_d$, and $\mlt(\h)$ is unitarily equivalent to $\M_d$.
\et
\bp
If the sequence $\log n_1, \log n_2, \ldots$ is linearly independent over $\mathbb Q$,
then the multiplier closure of $f(\bH_0)$ is $\B_d$ and $\h$ is weakly isomorphic as a reproducing kernel Hilbert space to $H^2_d$, by 
repeating  the proof of Theorem \ref{thm:universalPick}, replacing $p_i$ with $n_i$.

Conversely, assume that the sequence $\log n_1, \log n_2, \ldots$ is linearly dependent over $\mathbb Q$.
We show that the multiplier closure of
$f( \bH_0)$ is not $\B_d$ (and by the ideas of \cite[Section]{hartz15} this would also show that $\widehat{f(\bH_0)} \neq \B_d$).
We will exhibit a nonzero multiplier $q$ on $\B_d$ such that $q(f(s)) \equiv 0$.
Let $I$ and $J$ be disjoint finite subsets of the positive integers, and
let $\{\kappa_i\}_{i \in I \cup J}$ be nonnegative integers, not all zero, such that
\[
\sum_{i\in I} \kappa_i \log n_i  = \sum_{j\in J} \kappa_j \log n_j.
\]
Let $\mu$ be the multi-index supported on $I$ with $\kappa_i$ in the $i$th place, and let $\nu$ be defined likewise in terms of $J$.
Then we have that the polynomial $q(z) = b_\nu z^\mu - b_\mu z^\nu$ is not zero but satisfies
\begin{align*}
q(f(s)) &= b_\nu \Pi_{i\in I} b_{\kappa_i} n_i^{-s \kappa_i} - b_\mu \Pi_{j \in J} b_{\kappa_j} n_j^{-s \kappa_j}  \\
&= b_{\mu+\nu} (e^{-s\sum_{i\in I} \kappa_i \log n_i} - e^{-s\sum_{j\in J} \kappa_j \log n_j} ) \\
&= 0 .
\end{align*}
This shows that the multiplier closure of $f(\bH_0)$ is not all of $\B_d$.
\ep

\begin{remark}\emph{
When the sequence $\log n_1, \log n_2, \ldots, \log n_d$ is linearly independent over $\mathbb Q$, the space $\h$ is weakly isomorphic as a reproducing kernel Hilbert space to $H^2_d$, but these spaces are not isomorphic as reproducing kernel Hilbert spaces.
To illustrate the difference between these two notions, suppose that $d< \infty$, let $\h$ be as in Theorem \ref{thmm3}, and consider  $\varphi(s) = r - z_1 \circ f(s) = r - b_1 2^{-s}$ for $r \in (b_1, 1)$.
We have that $\varphi \in \mlt(\h)$ and $\inf_{s\in\bH_0}|\varphi(s)| \geq r - b_1 > 0$.
On the other hand, under the isomorphism of $\mlt(\h)$ and $\M_d$, $\varphi$ is mapped to $r - z_1$, which has a zero in $\B_d$, and is therefore not invertible in $\M_d$.
It follows that while $\inf_{s\in\bH_0}|\varphi(s)| \geq r - b_1 > 0$, there is no function $\psi \in \mlt(\h)$ such that $\varphi \psi =1$.
In contrast, if $\varphi : \B_d \to \C$ is a multiplier of $H^2_d$ and satisfies $\inf_{z\in\B_d}|\varphi(z)| > 0$ then $\varphi^{-1}$ is also a multiplier (this follows either from the corona theorem in Drury-Arveson space \cite{csw12}, or from the ``corona theorem for one function'' proved directly in \cite{fx13}).}
\end{remark}

Theorem \ref{thmm3} does not rule out that $\mlt(\h)$ is isometrically isomorphic to $\M_{d'}$ for $d'<d$, or for $d' = d$ in the case where $d = \infty$.
We now show that this possibility cannot happen.

\bt
If $\log n_1, \log n_2, \ldots$ are linearly dependent over $\mathbb Q$ then for any $d' \leq d$, $\mlt(\h)$ is  not isometrically isomorphic to $\M_{d'}$, and therefore $\h$ is not weakly isomorphic as a reproducing kernel Hilbert space to $H^2_{d'}$.
\et
\bp
Without loss of generality we assume that $d = \infty$.
Continue with the notation from the above theorem, and denote by $V$ the multiplier closure of $f(\bH_0)$ in $\B_\infty$.
To show that $\mlt(\h)$ is not isometrically isomorphic to $\M_{d'}$ for any $d' \leq d$, it suffices to show that $V$ is not an affine subspace of $\B_\infty$ (see \cite{salsha15}, Theorems 4.6 and 4.8).
Since $0 \in V$, it remains to show that $V$ is not a linear subspace.

Let $q$ be the polynomial defined in the proof of the previous theorem;
it is a nontrivial polynomial which vanishes on $V$, and
depends only on finitely many variables, say the first $N$ variables.
Let $\C^N \subset \ell^2$ be the finite dimensional subspace generated by the first $N$ standard basis vectors, and consider $q$ as a function on $\C^N$.
Since $q$ vanishes on $V$ and depends only on the first $N$ variables, it vanishes also on $P_{\C^N} V$.
The zero locus of $q$ in $\C^N$ has dimension $N-1$ (as a complex variety) and contains $P_{\C^N} V$, in particular $P_{\C^N} V$ is not equal to $\C^N$.
We will show that if $V$ is a subspace, then $P_{\C^N} V$ has dimension $N$ and hence is equal to $\C^N$ --- a contradiction.

Consider the sequence of functions $g_i : \bH_0 \to \C$ given by $g_i(s) = b_i e^{-s\log n_i}$.
Since $g_1, \ldots, g_N$ are linearly independent, the family of vectors
\[
\{(g_1(s), \ldots, g_N(s)) : s \in \bH_0\}
\]
spans all of $\C^N$.
If $V$ were a linear subspace, then $P_{\C^N} V$ would also be a linear space, therefore linear combinations of $P_{\C^N} f(\bH_0)$  would lie in $P_{\C^N} V$.
But $P_{\C^N} f(s) = (g_1(s), \ldots, g_N(s))$, so taking linear combinations we obtain that $P_{\C^N}V =  \C^N$.
This contradiction completes the proof.
\ep

Our results settle the problem of when $\h$ has a multiplier algebra which is isometrically isomorphic to $\M_d$ for some $d$.
The question of when the multiplier algebra is algebraically isomorphic (or boundedly isomorphic, which is the same due to semisimplicity --- see \cite[Lemma 51]{drs15}) remains open.

\begin{question}
Let $\h$ be as above, and suppose that $d'$ is the maximal number of multiplicatively independent integers in the sequence $n_1,n_2,\ldots, n_d$.
Is it true that $\mlt(\h)$ is isomorphic to $\M_{d'}$?
\end{question}

By the main results of \cite{apv,kmccs13} the answer is yes when $d<\infty$ and  $d'=1$ (see also \cite[Section 2.3.6]{ars08}, \cite[Section 6]{drs15} and \cite[Section 7]{dhs15}).

\bibliography{../../../references}

\end{document}